\def\draft{n}
\documentclass{amsart}
\usepackage{fullpage,amssymb,epic,eepic,epsfig,pb-diagram,lamsarrow,pb-lams}

\theoremstyle{plain}

\newtheorem{theorem}{Theorem}
\newtheorem{proposition}{Proposition}[section]
\newtheorem{lemma}[proposition]{Lemma}

\theoremstyle{definition}

\newtheorem{question}{Question}

\theoremstyle{remark}

\newtheorem{remark}[proposition]{Remark}

\def\printname#1{
	\if\draft y
		\smash{\makebox[0pt]{\hspace{-0.5in}
			\raisebox{8pt}{\tt\tiny #1}}}
	\fi
}

\newcommand{\psdraw}[2]
         {\begin{array}{c} \hspace{-1.3mm}
	\raisebox{-4pt}{\epsfig{figure=draws/#1.eps,width=#2}}
	\hspace{-1.9mm}\end{array}}

\newlength{\standardunitlength}
\catcode`\@=11
\long\def\@makecaption#1#2{%
     \vskip 10pt

\setbox\@tempboxa\hbox{
       \small\sf{\bfcaptionfont #1. }\ignorespaces #2}%
     \ifdim \wd\@tempboxa >\captionwidth {%
         \rightskip=\@captionmargin\leftskip=\@captionmargin
         \unhbox\@tempboxa\par}%
       \else
         \hbox to\hsize{\hfil\box\@tempboxa\hfil}%
     \fi}
\font\bfcaptionfont=cmssbx10 scaled \magstephalf
\newdimen\@captionmargin\@captionmargin=2\parindent
\newdimen\captionwidth\captionwidth=\hsize
\catcode`\@=12

\newcommand{\tr}{\operatorname{tr}}

\def\lbl#1{\label{#1}\printname{#1}}

\def\BZ{\mathbb Z}
\def\BQ{\mathbb Q}

\def\D{\Delta}

\def\R{\mathcal R}

\def\La{\Lambda}
\def\l{\lambda}

\def\S{\Sigma}

\def\la{\langle}
\def\ra{\rangle}
\def\i{^{-1}}

\def\w{\omega}

\def\i{^{-1}}

\def\a{\alpha}

\def\bd{\partial}
\def\g{\gamma}

\def\b{\beta}

\def\sub{\subset}

\renewcommand{\ker}{\operatorname{Ker}}

\def\iso{\cong}

\def\sminus{\smallsetminus}

\def\ti{\widetilde}
\def\Sei{\mathrm{Sei}}

\def\sarrow{\rightsquigarrow}
\def\Sone{\stackrel{S_1}\sarrow}
\def\Stwo{\stackrel{S_2}\sarrow}

\def\matt#1#2#3#4#5#6#7#8#9{\left(
\begin{matrix}
 #1 & #2 & #3 \\
 #4 & #5 & #6 \\
 #7 & #8 & #9   
\end{matrix}
\right)}

\def\ab{\mathrm{ab}}
\def\chiD{\chi_{\Delta}}
\def\ygraph{$\mathrm{Y}$-graph}

\def\Lhat{\hat\Lambda}

\def\Laab{\Lambda^{\mathrm{ab}}}
\def\Lhatab{\hat\Lambda^{\mathrm{ab}}}

\def\cyclic{\mathrm{cyclic}}
\def\qad{\mathcal R^{\mathrm{ad}}}

\begin{document}


\title[Analytic invariants of boundary links]{Analytic invariants of boundary 
links}

\author{Stavros Garoufalidis}
\address{School of Mathematics \\
          Georgia Institute of Technology \\
          Atlanta, GA 30332-0160, USA. }
\email{stavros@math.gatech.edu}
\author{Jerome Levine}
\address{Department of Mathematics\\
         Brandeis University\\
         Waltham, MA 02454-9110, USA. }
\email{levine@brandeis.edu}

\thanks{The  authors were partially supported by NSF grants
        DMS-98-00703 and DMS-99-71802 respectively, and by an Israel-US BSF
grant. 
        This and related preprints can also be obtained at
{\tt http://www.math.gatech.edu/$\sim$stavros } and
{\tt http://www.math.brandeis.edu/Faculty/jlevine/ }
\newline
1991 {\em Mathematics Classification.} Primary 57N10. Secondary 57M25.
}

\date{
This edition: October 19, 2000. \hspace{0.5cm} First edition: October 17, 
2000.}


\begin{abstract}
Using basic topology and linear algebra, we define a plethora of
invariants of boundary links whose values are power series with noncommuting
variables. These turn out to be useful and elementary reformulations of an
invariant originally defined by M.
Farber \cite{Fa2}.
\end{abstract}

\maketitle



\section{Introduction}
\lbl{sec.intro}
\subsection{History and Purpose.} 
\lbl{sub.history}

In a series of papers, M. Farber used homological methods to introduce 
an invariant of boundary links with values in a ring of rational functions 
with noncommuting variables \cite{Fa2}. 
A similar invariant to that of Farber was recently introduced by V. Retakh, 
C. Reutenauer and A. Vaintrob \cite{RRV} based on the notion of  
quasideterminants.

The purpose of this paper is to give an interpretation of Farber's
invariant as a simple invariant 
of the Seifert matrix of a boundary link, which is more elementary and makes
calculation  more straightforward. From this point of view we will, in fact,
define a whole spectrum of invariants  which take values in non-commutative 
power series rings. Although these invariants all turn out to be
determined by 
Farber's---see Theorem \ref{th.chi}---it is useful to have the
different formulations. An  example is given by $\chiD$---see page
\pageref{page.delta}---which has direct application to the study of the
Kontsevich integral of a boundary link and its rationality properties, as 
will be explained in subsequent publications \cite{Ga,GK}.  $\chiD$ also
gives a natural way to see that Farber's invariant determines the natural analog
of the Alexander polynomial for a boundary link (the classical Alexander
polynomial of a boundary link is $0$). See Proposition \ref{prop.2}.

We would like to thank Michael Farber for useful discussions.
\subsection{Boundary links and their refinements}
\lbl{sub.whatis}

All manifolds will be oriented and all maps will be smooth and orientation
preserving.
A {\em boundary link} ($\bd$-link) in a 3-manifold is an oriented link which is
the 
boundary
of a disjoint union of connected surfaces, each with one boundary component.
A choice of such surfaces is called a {\em Seifert surface}
of the boundary link. 
It is well-known that in the case of boundary links (unlike the case of knots)
the cobordism class, relative boundary, of a Seifert surface for a given link 
is not unique. There are at least two ways to overcome this difficulty,
as was explained by Cappell-Shaneson \cite{CS} and Ko \cite{K2}:
\begin{itemize}
\item 
A {\em $\S$-boundary link} $L$ (or simply, a {\em $\S$-link}) in a 3-manifold
$M$ is a choice,   up to isotopy, of Seifert surface $\Sigma$ in $S^3$
such that $\partial \Sigma=L$ .
\item 
An {\em $F$-boundary link} $L$ of $n$ components (or, simply, an 
{\em $F$-link}) is a link,   up to isotopy, equipped with a map $\phi:
\pi_1(M \sminus L)
\to F$ where $F$ is the free group on $n$ letters and $\phi$ maps a choice of 
meridians of $L$ to a   basis of $F$. $\phi$ is called a {\em splitting
map} for $L$.
\end{itemize}

It turns out that $F$-links can be identified with the set of cobordism 
classes, rel boundary, (or {\em tube equivalence} classes) of Seifert 
surfaces--- see Gutierrez and Smythe  \cite{Gu,Sm}.
Let $A_n$ denote the group of automorphisms $\a$ of the free group  
$F(t_1,\dots,t_n)$ that satisfy $\a(t_i)=w_i t_i w_i^{-1}$ for some $w_i \in
F(t_1,\dots,t_n)$, for all $i$, \cite{CS,K2}.  $A_n$
acts on the set of $F$-links by composition with the splitting map $\phi$. In
\cite{K2} a simple set of generators for $A_n$ was given, and the action of
these generators was described geometrically as what was there called {\em
cocooning}. It turns out that the set of 
equivalences classes of $F$-links, modulo the $A_n$ action, can be identified
with the set of $\partial$-links. 

We denote by $X_L^{\w}$ the $F$-covering of $S^3 -L$ associated with
$\ker\phi =\pi_{\w}$, the intersection of the lower central series of $\pi
=\pi_1 (S^3 -L)$.

\subsection{Seifert matrices of boundary links}
\lbl{sub.noncalculus}

There is an algebraic notion of a {\em Seifert matrix} associated to a 
$\S$-link of $n$ components, 
\cite{K1,K2}. These matrices are partitioned into $n\times n$ blocks
of matrices, corresponding to the link components. Let $\Sei(n)$ denote the set
of matrices $A=(A_{ij})$ of square matrix blocks $A_{ij}$ for $i,j=1\dots,n$, 
with integer entries, satisfying the conditions
$$
A_{ij}'=A_{ji} \,\, \text{for $i \neq j$ and} \,\,  \det(A_{ii}-A_{ii}')=1
\,\, \text{for all $i$}.
$$
Let $\Sei$ denote the set of all Seifert matrices. 
The Seifert matrix associated to a $\S$-link (resp. $F$-link, $\partial$-link)
is an element of $\Sei(n)$, well-defined up to $S_1$-equivalence
(resp. $S_{12}$-equivalence, $S_{123}$ equivalence), where
$S_1$ stands for {\em congruence}, $S_2$ stands for {\em stabilization}
and $S_3$ stands for equivalence under an algebraic action of $A_n$ on
$\Sei(n)$ defined by Ko \cite{K1,K2} (see Section \ref{sub.all} below). 
Note that $S_{123}=S_{12}$ for $n\le 2$, since $A_n$ consists entirely of inner
automorphisms which act trivially on Seifert matrices.
We have a commutative diagram 

$$
\divide\dgARROWLENGTH by2
\begin{diagram}
\node{\S-\text{links}}\arrow{e}\arrow{s}\node{\Sei(n)/(S_1)}\arrow{s} \\
\node{F-\text{links}}\arrow{e}\arrow{s}\node{\Sei(n)/(S_{12})}\arrow{s} \\
\node{\partial-\text{links}}\arrow{e}\node{\Sei(n)/(S_{123})} 
\end{diagram}
$$

Set $\La=\BQ[F]$ the group-ring with rational coefficients and  $\Lhat$ its
completion with respect to powers of the augmentation ideal. Then $A_L =H_1
(X_L^{\w},\BQ )$ is a $\La$-module. Let $\Laab
=\BQ[H]$, where $H$ is the free abelian group on generators $(t_1,\dots,t_n)$.
If $X_L^{\ab}$ denotes the universal abelian covering of $S^3 -L$, then $H_1
(X_L^{\ab},\BQ )$ is a $\Laab$-module. 
Note that $\Lhat$ can be identified
with the power series ring in the $n$ noncommuting variables $x_i =t_i -1$ and
$\Lhatab$ with the power series ring in $n$ commuting variables $x_i =t_i 
-1$. $\La$ (and also, $\Lhat, \Laab, \Lhatab$) are rings with (anti)-involution
given by $g\to \bar g=g^{-1}$ for $g \in F$. Note that $\bar x_i=-(x_i+1)^{-1}
x_i$. 
The action of $A_n$ on $F$ extends naturally to $\La$ and $\Lhat$ and 
induces the trivial action on $\Laab$.
Now, we can introduce analytic invariants of the set $\Sei$:
Let $f \in \BQ\la\la x,z \ra\ra$ be a noncommutative power
series in two variables. We will say $f$ is {\em admissible} if, for any
non-negative integer $n$, there are only a finite number of terms in $f$ of 
total $x$-degree $n$. The admissible power series form a subring $\qad$ of 
$\BQ\la\la x,z
\ra\ra$. Now let $X=\text{diag}(x_1,\dots,x_n)$ be a (block) diagonal
matrix. Then we let 
 $$\chi_f: \Sei(n) \to \Lhat
\text{ be defined by } \,\,\, \chi_f(A)=\tr(f(X,Z_A)-f(X,I_{1/2}))
$$
where $Z_A=A(A-A')^{-1}$ and $I_{1/2}$ is the block diagonal matrix in which
  half of the diagonal entries  in each diagonal block
  are $0$ and half are $1$. Note that
$f(X,I_{1/2}))$ is independent of how the $0$'s and $1$'s are distributed

\begin{theorem}
\lbl{thm.all}
For all admissible $f$, $\chi_f$ descends to a map
$$
\Sei(n)/(S_{12})\to \Lhat.
$$
\end{theorem}


\begin{remark} If $f\in\BZ\la\la x,z\ra\ra$ then $\chi_f (A)$ has integer
coefficients.
\end{remark}
\begin{question}
If $\a\in A_n$, $f\in\qad$ and $A$  a Seifert matrix, is $\chi_f (\a\cdot A)$
determined by $\chi_f (A)$ and $\a$?
\end{question}

Let $\R(a_1,\dots,a_n)$ denote the   subring of $\BQ\la\la a_1,\dots,a_n\ra\ra$
consisting of the {\em rational functions} in the 
noncommuting variables $\{a_1,\dots,a_n\}$ (see \cite{B}). This can be defined
as  the smallest subring of $\BQ\la\la a_1,\dots,a_n\ra\ra$ containing the
polynomials $\BQ [a_1,\dots,a_n ]$ and closed under the operation of taking
inverses of {\em special} series, i.e. those  $f$ with constant term $f(0,\dots
,0)=1$\cite[p.6]{B}. Let
$\R_*  (x,z)$ denote the smallest subring of $\R(x,z)$ containing the
polynomials and closed under the operation of taking inverses of {\em
extra-special} series, i.e. admissible $f(x,z)$ which satisfy $f(0,z)=1$.
Clearly $\R_* (x,z)\sub\R (x,z)$. We also note that $\R_* (x,z)\sub\qad$ since
it is not hard to see that if $f$ is special then $f\i$ is admissible if and
only  $f$ is  extra-special.

\begin{question}
Is $\R_* (x,z)=\R (x,z)\cap\qad$?
\end{question}

\begin{proposition}\lbl{prop.rat}
If $f(x,z)\in\R_* (x,z)$, then $\chi_f\in\R (x_1 ,\dots ,x_n )$.
\end{proposition}

$\chi_f$ satisfies a general duality property. 
Define  anti-involutions $f\to\ti
f$  and $f\to\bar f$ on $\BQ\la\la x,z \ra\ra$ and $\BQ\la\la
x_1 ,\dots ,x_n \ra\ra$ by the properties 
$\ti x=x, \bar x=-x(1+x)\i , \ti z=\bar z =z, \ti x_i =x_i ,\bar x_i =-x_i
(1+x_i )\i $ 
and $\widetilde{fg}=\ti g\ti f ,\overline{fg}=\bar g\bar f$ , and an involution 
$f\to\hat f$ on $\BQ\la\la
x,z \ra\ra$ and $\BQ\la\la x_1 ,\dots ,x_n \ra\ra$,  by $\hat x=\bar
x, \hat z=z ,\hat x_i =\bar x_i$ and $\widehat{fg}=\hat f\hat g$. Note that the
composition of any two of the maps $f \to \tilde{f}, \bar f$ or $\hat{f}$
on $\BQ\la\la x_1 ,\dots ,x_n \ra\ra$ is equal to the third.

\begin{proposition}\lbl{prop.dual}
For admissible $f$, we have that 
$\ti\chi_{f(x,z)}=\chi_{\ti f(x,1-z)}$ and $\chi_{\hat f}=\hat\chi_f$. 
Therefore we also have $\bar\chi_{f(x,z)}=\chi_{\bar f(x,1-z)}$.
\end{proposition}

Note that if $f$ is admissible then
\begin{itemize}
\item $\ti f, \hat f$ are admissible, and
\item  $f(x,1-z)$ is
defined (which is not true for every $f\in\BQ\la\la x,z \ra\ra$) and 
admissible.
\end{itemize}

Let $f_{\D}=\log\left(xz+1\right) 
\in \BQ\la\la x,z\ra\ra$. Note that $f_{\D}$ is admissible---in fact any
$f\in \BQ\la\la x,z\ra\ra$ of the form $f(x,z)=G(xz)$, where $G(y)\in\BQ\la\la
y\ra\ra$, is admissible. Let us denote \lbl{page.delta} $\chi_{f_{\D}}$ by
$\chiD$.
Our interest in $\chiD$ comes from the fact that it can be identified with
the ``wheels part'' of a (version of) the Kontsevich integral of
$F$-links, as will be explained in a separate publication, \cite{GK}.
For now, let us explain the relation between $\chiD$ and the algebraic
topology of the complement of a boundary link $L$. 

Let  $\D^b(L) \in \Laab/(\mathrm{units})$, denote the order of the torsion 
$\Laab$-module  
$A^{\ab}_L=\text{torsion}_{\Laab}H_1(X_L^{\ab}, \BQ)$.

It is well-known that $\D^b =\D^b
(L)$ satisfies:
\begin{enumerate}
\item $\D^b (1,\dots,1)=\pm 1$ and 
\item $\D^b (t\i _1 ,\dots,t\i _n )$ is a
unit multiple of $\D^b (t_1,\dots,t_n )$ in $\Laab$.
\end{enumerate}
It follows that we can choose a unique
(normalized)  representative in $\Laab$ such that 
\begin{enumerate}
\item[(1')] $\D^b (1,\dots,1)=1$ and 
\item[(2')] $\D^b (t\i _1 ,\dots,t\i _n )=\D^b (t_1,\dots,t_n )$.
\end{enumerate}
We call this normalized representative the {\em torsion polynomial} of $L$.
\begin{proposition}
\lbl{prop.2}
\begin{itemize}
\item[(a)] (Abelianization)   If $\chiD^{\ab}$ denotes the abelianization of
$\chiD$,
then
$$
\chiD^{\ab}=\log \D^b \in \Lhatab. 
$$
where $\D^b$ is the torsion polynomial.
\item[(b)] (Realization) For every element $\l\in\La$ with integer
coefficients  satisfying $\l
(1,\dots,1)=1$ and $\l=\bar \l$, there exists an $F$-link $L$ with 
$H_1 (X^{\omega}_L ,\BZ)\iso\La /(\l )$, where $(\l)$ denotes the left ideal 
generated by $\l$. As a consequence every element $\D$ of $\Laab$ satisfying 
(1') and (2') can be realized as the torsion polynomial of some boundary link. 
\item[(c)] (Duality)  
  $ \chiD= \bar{\chi}_{\Delta}$ in $\Lhat/(\cyclic)$, 
the quotient of
$\Lhat$ by its {\em subgroup} generated by $(ab-ba)$, for $a,b \in \Lhat$. 
\end{itemize}
\end{proposition}

Thus, $\chiD^{\ab}$ determines the  torsion
polynomial. In contrast, the classical multivariable Alexander polynomial
of a boundary link vanishes, and in general it is not known which Laurent
polynomials can be realized as the multivariable Alexander polynomials of a
link.

For an $F$-link $L$, we can think of $\chiD$ as an analogue of the order
of the $\La$-module $A_L$ (even though the notion of order does not make sense
for  
$\La$-modules). 

\begin{proposition}
\lbl{prop.compare}
If $\Phi (x,z)=(xz+1)^{-1}x \in \BQ \la \la x,z \ra \ra$, then $\chi_{\Phi}$
is related to 
Farber's $\chi$-function \cite[Section 2.4]{Fa2} by the formula
$$\chi -\chi_{\Phi}=\sum_{i=1}^n g_i \left(x_i - \bar x_i \right)$$
for some 
non-negative 
integers $g_i$.
\end{proposition}

\begin{remark} 
It follows from Farber's approach that $\chi$  only depends on the $\La$-module
$A_L$, and, therefore, this is also true for $\chi_{\Phi}$, when
$\Phi(x,z)=(xz+1)^{-1}x$, 
since the integers $g_i$ are half the ranks of the
$x_i$-components of the {\em minimal lattice} in $A_L$, as is
demonstrated in the proof of the above proposition. 
\end{remark}

\begin{question}
\lbl{que.blanchfield}
Is there some way to see directly that $\chi_{\Phi}$ depends only on $A_L$?
\end{question}

In \cite{Fa2} it is shown that $\chi$ (and thus $\chi_{\Phi}$)
satisfies the duality property $\chi +\bar\chi =0$.  
We reprove this using Proposition \ref{prop.compare} and \ref{prop.dual}.

\begin{proposition}
\lbl{prop.far}
  For any $F$-link, we have 
$$\chi_{\Phi}=-\bar\chi_{\Phi}$$
\end{proposition}

\begin{question} (Realization)  Can every {\em rational} power
series $\rho$,  with integer coefficients, satisfying the duality property
$\rho =-\bar\rho$
 be realized as $\chi_\Phi(L)$
for some $F$-link $L$? 
\end{question}
\begin{question} For the cyclic module in Proposition \ref{prop.2}(c), what is
$\chi$?
\end{question}

It is interesting, if perhaps disappointing, that this array of invariants
are actually all determined by the original $\chi$ of Farber.
\begin{theorem}\lbl{th.chi}
For any $f\in\BQ\la\la x,z\ra\ra, \chi_f (L)$ is completely determined by $\chi
(L)$, and therefore depends only on $A_L$.
\end{theorem}

\begin{remark}
It is pointed out in \cite[Prop. 5.2]{Fa2} that $\chi (L)$ determines $A_L$
when it is semi-simple but not otherwise. For example, it follows from
\cite[Prop. 2.5(c)]{Fa2} that $\chi$ is not sensitive to different extensions 
of the same modules. In particular, for a knot $K$, $\chi (K)$ is determined 
by the Alexander polynomial \cite[Section 10.4]{Fa2} and it is well-known 
that there exist knots with the same Alexander polynomial but different 
Alexander modules.
\end{remark}


Finally we consider some examples. If $L$ is an $F$-link, let $L'$ denote the
reflection (sometimes called {\em mirror image}) of $L$ with the natural 
$F$-structure induced from that of $L$ by 
the automorphism of $F$ defined by $t_i\to t_i\i$. If $A$ is a Seifert matrix 
of $L$ then $A'$ is a Seifert matrix for $L'$. Note that 
$Z_{A'}=I-Z_A =SZ'_A S\i$.

\begin{proposition}\lbl{prop.reflect}
 For any $f\in\BQ\la\la x,z\ra\ra$, $\chi_f
(L')=\ti\chi_{\ti f}(L)$. In particular $\chi (L')=\ti\chi (L)$. 
\end{proposition}
For $2$-component links we have not been able to find any examples such that
$\chi (L')\not=\chi (L)$.
\begin{question}
Is $\ti\chi (L)=\chi (L)$ for any $2$-component $F$-link?
\end{question}
On the other hand for $3$-component $F$-links it is not hard to find such 
examples.

\begin{proposition}\lbl{prop.refl}
There exist $3$-component $F$-links such that $\chi(L')=\ti\chi(L)\not=
\chi (L)$.
\end{proposition}

\section{Proofs}
\lbl{sec.proofs}

\subsection{Proof of Theorem \ref{thm.all}}
\lbl{sub.all}

Let us introduce three moves on the set $\Sei$ of Seifert matrices:
\begin{itemize}
\item[$S_1:$]
Replace $A$ by $PAP'$ for a block diagonal matrix $P=\text{diag}
(P_1,\dots,P_n)$ of unimodular matrices $P_i$ with integer entries.
\item[$S_2:$]
Replace $A$ by 
$$
\left(
\begin{matrix}
A & \rho & 0 \\
\rho'  & 0    & 1 \\
0 & 0    & 0
\end{matrix}
\right)
\,\,\, \text{ or }
\left(
\begin{matrix}
A & \rho    & 0 \\
\rho' & 0    & 0 \\
0 & 1    & 0
\end{matrix}
\right)
$$
for a  column vector $\rho$, where, for some $i$,  the two new rows are
added to $A_{ij}, 1\le j\le n$ and the two new columns are added to $A_{ji},
1\le j\le n$.
\item[$S_3:$]
The move that generates $A_n$-equivalence, where the algebraic action of
$A_n$ on $\Sei(n)$ is described in \cite{K1,K2}.
\end{itemize}
Note that $S_1,S_2$ generate the so-called $S$-equivalence of Seifert
matrices. 

  Given a Seifert matrix $A$, we define
$Z_A=A(A-A')^{-1}$ and $S_A=A-A'$
(or simply, $Z$ and $S$ in case $A$ is clear), following Seifert.
Note that $S$ is block-diagonal. The behavior of $Z$ under $S$-equivalence of 
$A$ is described by the following elementary matrix calculation

\begin{lemma}
\lbl{lem.ZS}
If $A \Sone B=PAP'$, 
then $Z_B= PZ_AP^{-1}$. \newline
If $A \Stwo B$, then 
$$
Z_B=
\left(
\begin{matrix}
Z_A & 0 & \star \\
\star & 1    & \star \\
0 & 0    & 0
\end{matrix}
\right)
\,\,\, \mathrm{ or } \,\,\,
\left(
\begin{matrix}
Z_A & \star    & 0 \\
0 & 0    & 0 \\
\star & \star    & 1
\end{matrix}
\right).
$$
\end{lemma}

\begin{proof}[Proof of Theorem \ref{thm.all}]
We need to show that $\chi_f$ is invariant under the moves $S_1$ and
$S_2$.
If $A \Sone B$, then $f(X,Z_B)=
P f(X,Z_A) P^{-1}$ thus $\chi_f(B)=\chi_f(A)$. 
If $A \Stwo B$, then the following identity 
$$
\left(
\begin{matrix}
C  & 0  & \star  \\
\star  & c & \star   \\
0  &    0 &  0 
\end{matrix}
\right)
\left(
\begin{matrix}
C'  & 0  & \star  \\
\star  & c' & \star   \\
0  &    0 &  0 
\end{matrix}
\right)
=
\left(
\begin{matrix}
C C' & 0  & \star  \\
\star  & c c' & \star   \\
0  &    0 &  0 
\end{matrix}
\right)
$$
implies that $\chi_f(B)=\chi_f(A)$. 
\end{proof}

\begin{lemma}
\lbl{lem.presentation}
Given an $F$-link $L$ with a Seifert matrix $A$, then  $XZ_A+I$ is a
presentation matrix for $A_L$ over $\La$, and for  
$A_L^{\ab}$ over $\La^{ab}$. 
\end{lemma}

\begin{proof}
It is well-known (see \cite{K2}) that a presentation matrix for 
$A_L$ is $TA-A'$, and similarly for $A_L^{\ab}$. Since
$TA-A'=(T-I)A+ (A-A')=(XZ+I)(A-A')$, the lemma follows. 
\end{proof}

\subsection{Proof of Proposition \ref{prop.rat}}
Let $\R '(x,z)$ denote the subring of $\qad$ consisting of all admissible $f$
such that, for any scalar matrix $Z$ of the appropriate size, $f(X,Z)$ is a
matrix all of whose entries are rational in $\BQ\la\la x_1,\dots ,x_n\ra\ra$. 
It will suffice to show that $\R_* (x,z)\sub\R '(x,z)$. To prove this we need 
to show that if  $f\in\R '(x,z)$ is extra-special, then $g=f\i\in\R'(x,z)$. 

Consider the matrix equation $f(X,Z)g(X,Z)=I$. This defines a system of
equations
for the entries of $g(X,Z)$ of the form
$$\sum_j a_{rj}y_j =b_r $$
where $a_{rj}, b_r\in\R (x_1 ,\dots ,x_n )$. Since $f$ is extra-special,
$f(0,Z)=I$, which implies that, with the correct choice of numbering of the
equations, $a_{rj}(0,\dots
,0)=\delta_{rj}$. Now we can apply \cite[Proposition 2.1]{FV} to
conclude that the solutions $y_r$, which are the entries of $g(X,Z)$, are unique
and rational.

\subsection{Proof of Proposition \ref{prop.dual}} 
 It follows from the definition of $Z$ that $Z+SZ'S^{-1}=I$, 
where $S=A-A'$. This implies that
\begin{equation}\lbl{eq.dual}
Z=S(I-Z')S^{-1}
\end{equation}
Thus 
$$\tr f(X,Z)=\tr f(X,I-Z')=\ti\tr\ti f(X,1-Z) $$
using the facts that $S$ commutes with   $X$, that $\tr Y=\tr Y'$, for any 
square matrix $Y$ and that $\tr (WY)=\tr (YW)$, if the entries of $W$ commute
with
the entries of $Y$. From this we deduce the first equality.

The second equality is clear.

\begin{remark}
\lbl{rem.dual}
Let $[\qad,\qad]$ denote the abelian subgroup of the ring $\qad$ generated
by $fg-gf$ for $f,g \in \qad$. It is easy to see that for all $f \in 
[\qad,\qad]$ we have $\chi_f=0 \in \Lhat/(\cyclic)$.
\end{remark}

\subsection{Proof of Proposition \ref{prop.2}}

To prove (a) first note that the normalized $\Delta^b$ can be defined by the
equation 
$$\Delta^b =\det (T^{1/2}A-T^{-1/2}A')=\det ((I+X)^{-1/2}(I+XZ))$$
Thus we have 
\begin{eqnarray*}
\log\Delta^b &=&\tr\log ((I+X)^{-1/2}(I+XZ))\\
&=&\tr\log (I+X)^{-1/2} +\tr\log (I+XZ)\\
&=&\tr\log (I+XZ) -\tfrac{1}{2}\tr\log (I+X)=\chi_{\Delta}
\end{eqnarray*}
This uses the following lemma, which is probably well-known.
\begin{lemma} Suppose that $M, M_1 ,M_2$ are matrices of the form $I+N$ over a
completed commutative power series ring, where $N$ has all entries of degree
$>0$. Then we have the following identities.
\begin{eqnarray}
\tr\log (M_1 M_2 )&=&\tr\log (M_1 )+\tr\log (M_2)\lbl{eq.mult}\\ 
\log \det(M)&=&\tr\log (M) \lbl{eq.detr}
\end{eqnarray}
\end{lemma}
\begin{proof}\eqref{eq.mult} follows from the Campbell-Baker-Haussdorf formula
and the fact that $\tr (AB)=\tr (BA)$ if $A,B$ are matrices over a commutative
ring.

To prove \eqref{eq.detr}, first note that it is obvious if $M$ is triangular.
Secondly, it follows from \eqref{eq.mult} that if it is true for $M_1$ and
$M_2$, then it is true for $M_1 M_2$. Thus if will follow from the fact that any
such $M$ can be written $M=LU$, where $L$ is lower triangular (i.e. $l_{ij}=0$
if $i<j$) and $U$ is upper triangular. We prove this by induction on the size of
$M$. 

Write $M=\left(\begin{array}{rr} u & \b \\\a & \ti M\end{array}\right)$, where
$\a$
is a column vector and $\b$ is a row vector. By induction we can write $\ti
M-u\i\a\b =\ti L\ti U$, for triangular matrices $\ti L,\ti U$. Now we define
$$L=\left(\begin{array}{rr} u & 0\\\a & \ti L \end{array}\right)
\quad\text{and}\quad
U=\left(\begin{matrix}1 & u\i\b\\0 & \ti U \end{matrix}\right)$$
One checks immediately that $M=LU$.
\end{proof}
(b) follows from a general construction in \cite{Le}. Consider the trivial link
$L_0\sub
S^3$ with $n$ components. Then the splitting map $\phi$ is an isomorphism.
Consider the universal cover $X_{L_0}^{\w}$ of $S^3 -L_0$. Given 
$\l =\sum_{g\in F}a_g
g$ satisfying $\l =\bar\l$, we can construct a simple closed curve $\g$ in $S^3
-L_0$ which is null-homotopic and unknotted in $S^3$ such that, if $\ti\g$ is
any lift of $\g$ in $X_{L_0}^{\w}$ then the linking numbers of $\ti\g$ and its
translates is given by 
$$ \text{lk}(\ti\g , g\ti\g )=a_g \text{\ if }g\not= 1 $$
This construction is described in \cite{Le}. Now do a $+1$-surgery on $\g$ to
produce $\S^3$, which, since $\g$ was unknotted, is diffeomorphic to $S^3$. Let
$L\sub S^3$ be the link corresponding to $L_0\sub\S^3$ under such a
diffeomorphism. Note that surgery on all the lifts of $\g$ produces an 
$F$-covering of $S^3 -L$ and so $L$ is canonically a $F$-link. The argument 
in \cite{Le} shows that $H_1 (X_{L}^{\w})\iso\La /(\l)$.

For (c), 
we will use  Proposition \ref{prop.dual} and Remark \ref{rem.dual}.
Since $f_\D(x,z)=\log(1+xz)$, it is easy to see that   
 $f_\D(x,z)=f_\D(z,x) \bmod [\qad,\qad]$. 
On the other hand, we have  
\begin{eqnarray*}
\ti f_\D(x,1-z) &=& \log(1+(1-z)x)=\log(1+x-zx)=\log((1+z\bar x)(1+x)) \\
&=&
\log(1+z\bar x)+\log(1+x) \bmod [\qad,\qad]  \\ &=&
\bar f_\D(x,z)+\log(1+x)  \bmod [\qad,\qad]  \\
&=&
\log(1+\bar x z)+\log(1+x) \bmod [\qad,\qad] \\
&=&
\hat f_\D (x,z)+\log(1+x) \bmod [\qad,\qad].
\end{eqnarray*}
Proposition \ref{prop.dual} and Remark \ref{rem.dual} imply that
$$\chi_{f_\D(x,z)}=\ti\chi_{\ti f_\D(x,1-z)}=\ti\chi_{\hat f_\D (x,z)}
+\chi_{\log(1+x)}=\bar\chi_{f_\D}+\chi_{\log(1+x)}
\in \Lhat/(\cyclic).
$$
Since $\chiD(x_i)=\chi_{f_\D(x,z)}$  and $\chi_{\log(1+x)}=0$, it follows
that $\chiD=
\bar\chi_\D
\in \Lhat/(\cyclic)$.
\qed

\subsection{Proof of Proposition \ref{prop.compare}}

 Let $A$ be any Seifert matrix for $L$. We can construct a higher-dimensional
{\em simple} link $\ti L$ in $S^{4k+3}$, for some large $k$, which has a 
Seifert manifold $W$ yielding $A$ as its Seifert matrix (see, e.g. \cite{K2}).
Since the Seifert matrix determines
the link module, via the presentation matrix $TA-A'$ we have
$A_L=H_1(X^\omega_L,\BQ)\iso A_{\ti L}= H_{2k+1} (X^\omega_{\ti L},\BQ)$.
Therefore
$\chi$ for $A_L$ is the same as $\chi$ for $A_{\ti L}$. Now we can do surgery on
$W$
to obtain a {\em minimal Seifert manifold} $V$ 
for $\ti L$, whose components are
$2k$-connected---see  \cite[Section 6.12]{Fa1} and \cite{Gu}. This
determines a minimal
lattice $J$ for $A_{\ti L}$, according to \cite[p.563-4]{Fa2}. The Seifert
matrix
$B$ determined by $V$ is S-equivalent to $A$ and so $\chi_{\Phi}(L)= \tr
((I+XZ)\i X)- \tr ((I+XI_{1/2})\i X)$, where $Z=B(B-B')\i$. Note that
$J=\bigoplus_i x_i J$ and each $x_i J$ is isomorphic to $H_{2q+1}(V_i )$, where
$V_i$ is the
$i$-th component of $V$. Then, if $2g_i =\text{rank }H_{2q+1}(V_i )$, it is
straightforward to check that $ \tr ((I+XI_{1/2})\i X)=\sum_i g_i
(x_i -\bar x_i )$. The proof will be completed if we show that $\chi = \tr
((I+XZ)\i X)$.

Now $A_L$ is the $\La$-module with presentation matrix $XZ+I$,
as in Lemma \ref{lem.presentation}. The generators $\a_r$ of $A_L$,
corresponding to the columns of $XZ+I$, span the minimal lattice $J$ as
described in \cite[p.564]{Fa2}, since $B$ comes from a minimal Seifert manifold.
If we let $M_i=x_i A_L$, then the generators corresponding to the $i$th column
block of $XZ+I$ generate
$M_i$ since, if $\a_r$ denotes a generator corresponding to a column in the
$i$th column block, the $r$-th row of $XZ+I$ gives the relation 
$\a_r=-x_i \sum Z_{rs}\a_s$. Thus, $\pi_i$ is given by the matrix
$P_i=\matt 0 {} {} {} I {} {} {} 0$
where $I$ is in the $(i,i)$ block, $z$ is given by the matrix $Z'$ and 
$\partial_i$ is given by the matrix whose $i$th column block is the $i$th column
block of 
$-Z'$ and the other columns are zero---call this matrix $Z_i$.
Now, $\chi$ is given by
$$
\chi =\sum_k \sum_n \tr(\pi_k \partial_{a_1} \cdots \partial_{a_n})
x_{a_n} \cdots x_{a_1} x_k.
$$
But $\pi_k \partial_{a_1} \cdots \partial_{a_n}$ is given by the matrix 
$P_k Z_{a_1}' \cdots Z_{a_n}'$. Note that 
$Z_{a_1}'\cdots Z_{a_n}'$ is the matrix with only the $a_n$-th column block
nonzero and
the $(r,a_n)$ block is $(-1)^n Z_{r,a_1}' Z_{a_1,a_2}' \cdots 
Z_{a_{n-1},a_n}'$.
Multiply this by $P_k$, giving a matrix whose only nonzero entries are in the
$(k,a_n)$ block, and equal to
$ (-1)^n Z_{k,a_1}' Z_{a_1,a_2}' \cdots Z_{a_{n-1},a_n}'$.
Thus, we have a nonzero trace only if $k=a_n$, giving
\begin{eqnarray*}
\chi &=&
\sum_n \sum_{a_1,\dots,a_n} 
(-1)^n \tr( Z_{a_n,a_1}' Z_{a_1,a_2}' \cdots Z_{a_{n-1},a_n}')
x_{a_n} \cdots x_{a_1} x_{a_n} \\
&=&
\sum_n \sum_{a_1,\dots,a_n} 
(-1)^n \tr( Z_{a_n,a_{n-1}} Z_{a_{n-1},a_{n-2}} \cdots Z_{a_1,a_n})
x_{a_n} \cdots x_{a_1} x_{a_n} \\
&=&
\sum_n \sum_{a_1,\dots,a_n,a'_n} 
(-1)^n \tr( (XZ)_{a_n,a_{n-1}} (XZ)_{a_{n-1},a_{n-2}} \cdots (XZ)_{a_1,a'_n}
X_{a'_n,a_n}) \\
&=&
\sum_n (-1)^n \tr((XZ)^nX)=\tr\left((XZ+I)^{-1}X\right). 
\end{eqnarray*}

\qed

\subsection{Proof of Proposition \ref{prop.far}}

It is easy to see that $\ti\Phi(x,z)=\Phi (x,z)$.
Furthermore, since $\tr\Phi(X,I_{1/2})$ satisfies the asserted duality
statements, we can omit this part of the definition of $\chi_{\Phi}$ in the
following.
We have:
\begin{eqnarray*}
\ti\Phi(x,1-z) &=& (1+x(1-z))^{-1}x=(1+x-xz)^{-1}x=(1+\bar
xz)^{-1}(1+x)^{-1}x \\ &=&
-(1+\bar xz)^{-1}\bar x 
= -\Phi (\bar x,z)=-\hat\Phi (x,z).
\end{eqnarray*}
Proposition \ref{prop.dual} implies that  
$\chi_{\Phi}=\ti\chi_{\ti\Phi (x,1-z)}=
-\ti\chi_{\hat\Phi}=-\ti{\hat\chi}_{\Phi}=-\bar\chi_{\Phi}$.
\qed

\subsection{Proof of Theorem \ref{th.chi}}

 It suffices to consider the case where $f$ is a monomial, say 
$$f=x^{f_0}z^{e_1}x^{f_1}\cdots z^{e_k}x^{f_k} $$
where $e_i >0$ for $1\le i\le k$ and $f_i >0$ if $0<i<k$.  Note that we have a
general formula
\begin{equation}\lbl{eq.trace}
\tr f(X,Z)=\sum_{i_1 ,\dots ,i_k}\tr (Z^{e_1})_{i_1 i_2}(Z^{e_2})_{i_2
i_3}\cdots (Z^{e_k})_{i_k
i_1}x^{f_0}_{i_1}x^{f_1}_{i_2}\cdots x^{f_{k-1}}_{i_k}x^{f_k}_{i_1} 
\end{equation}
where $(Z^e )_{ij}$ denote the $(i,j)$-block of $Z^e$. Now we associate with 
$f$ another monomial $f'\in\BQ\la\la x,y,z\ra\ra$ by replacing each $z^{e_i}$ 
in $f$ by $(zy)^{e_i -1}z$, for every $1\le i\le k$ and replacing each 
$x^{f_i}$ by $x$, for every $0\le i\le k$ (even when $f_0$ or $f_k$ is zero). 
Now consider $\tr f'(X,Y,Z)$, where $Y=\text{diag}(y_1 ,\dots ,y_n )$ is a 
matrix identical to $X$
in which each $x_i$ is replaced by a new variable $y_i$. It is not hard to see,
using equation \eqref{eq.trace}, that  $\tr f'(X,Y,Z)$ and $f$ determine $\tr
f(X,Z)$ by replacing each $x_j$ in $\tr f'(X,Y,Z)$ with the appropriate power 
of $x_j$ and each $y_j$ by $1$. Furthermore, again using equation 
\eqref{eq.trace}, $\tr f'(X,X,Z)$ and $f$ determine $\tr f'(X,Y,Z)$ since $f$ 
tells us which $x_i$ in $f'(X,X,Z)$ to replace by $y_i$ to obtain 
$\tr f'(X,Y,Z)$. 

Finally we note that $f'(x,x,z)$ is a monomial of the form $(xz)^k x$, and so
coincides, up to sign, with the degree $k+1$ part of $\Phi$. Thus $\tr
f'(X,X,Z)$ is determined by $\tr \Phi (X,Z)$. This completes the proof.
\qed


\subsection{Proof of Propositions \ref{prop.reflect} and \ref{prop.refl}}
Since 
$\tr f(X,SZ'S\i )=\tr f(X,Z')=\ti\tr{\ti f }(X,Z)$,
 we conclude that $\chi_f (L')=\ti\chi_{\ti f}(L)$.
Since $\ti\Phi =\Phi$, it follows from Proposition \ref{prop.compare} that
$\chi(L')=\ti\chi (L)$. This proves Proposition \ref{prop.reflect}.

For Proposition \ref{prop.refl} let us consider the matrix
$$
A=\left(\begin{array}{rrr} M &-S &-S\\ S & M & -S \\ S & S &
M\end{array}\right) 
$$
where $S=\left(\begin{array}{rr} 0 & 1\\-1 & 0\end{array}\right)$ and $M$ is 
any $2\times 2$-matrix satisfying $M-M'=S$. This is a Seifert matrix for some 
$F$-link $L$. 

Then 
$$ Z_A =\left(\begin{array}{rrr} N&-I&-I\\I&N&-I\\I&I&N\end{array}\right) $$
where $N=MS$.
From the general formula
$$\tr_{\Phi}(Z)=\sum\tr x_{i_1}Z_{i_1 i_2}x_{i_2}\cdots
x_{i_k}Z_{i_{k}i_1}x_{i_1}=\sum (\tr Z_{i_1 i_2}\cdots  
Z_{i_{k}i_1})x_{i_1}x_{i_2}\cdots x_{i_k}x_{i_1}$$
which follows from equation \eqref{eq.trace}, 
we can see that $\tr_{\Phi}(Z_A )=\sum_{\bf m}a_{\bf m}\bf m$, summing over
non-commutative monomials ${\bf m}=x_{i_1}x_{i_2}\cdots x_{i_k}x_{i_1}$, where 
$a_{\bf m}=(-1)^r\tr M^s$ with $r=\#\{ j:i_{j+1}>i_j\}$ and $s=\#\{ j:i_j
=i_{j+1}\}$. Thus, for example if ${\bf m}=x_1 x_2 x_3 x_1$ then $a_{\bf m}=2$
whereas if ${\bf m}=x_1 x_3 x_2 x_1$ then $a_{\bf m}=-2$. Thus
$\chi_{\Phi}(A)\not=\ti\chi_{\Phi}(A)$ and we have the desired example.

\ifx\undefined\bysame
	\newcommand{\bysame}{\leavevmode\hbox
to3em{\hrulefill}\,}
\fi

\end{document}